\documentstyle[11pt]{article}
\begin{document}

\begin{center}{\bf  ON THE (1+3) THREADING OF SPACETIME WITH RESPECT TO AN ARBITRARY TIMELIKE VECTOR FIELD} \end{center}
\begin{center}{AUREL BEJANCU \\ Department of Mathematics \\ Kuwait University \\ P.O.Box 5969,
Safat 13060 Kuwait \\ E-mail:aurel.bejancu@ku.edu.kw \vspace{6mm}
\\ CONSTANTIN C\u ALIN \\ Department of Mathematics \\
Technical University "Gh.Asachi"\\ B-dul Mangeron no.67, 700050, Romania \\ E-mail:c0nstc@yahoo.com }\end{center}

\begin{abstract}
{We develop a new approach on the (1+3) threading of spacetime $(M, g)$ with respect to a congruence of curves
defined by an arbitrary timelike vector field. The study is based on spatial tensor fields and on the Riemannian
spatial connection $\nabla^{\star}$, which behave as $3D$ geometric objects. We obtain new formulas for local
components of the Ricci tensor field of $(M, g)$ with respect to the threading frame field, in terms of the
Ricci tensor field of $\nabla^{\star}$ and of kinematic quantities. Also, new expressions
for time covariant derivatives of kinematic quantities are stated. In particular, a new form of
Raychaudhuri's equation enables us to prove Lemma 6.2, which completes a well known lemma used in the proof of
Penrose-Hawking singularity theorems.Finally, we apply the new $(1+3)$ formalism to the study of the dynamics
of a Kerr-Newman black hole.}
 \end{abstract}


%\noindent {\bf PACS: 04.20.-q, 04.20.Cv, 04.70.-s}

\medskip
\newpage
\section{Introduction}
The (1+3) threading  of spacetime by a congruence of curves determined by a unit timelike vector field $\xi$
(4-velocity) is by now a well established theory which studies the geometry, dynamics, and observational properties
of some well known cosmological models. Most of the important results on this theory, and an exhaustive list of
 references wherein these results have been published, can be found in the excellent monograph of Ellis, Maartens
and MacCallum \cite{ej}.\par
Our work on this matter is motivated by the simple remark that it is difficult to apply the above theory to
the metrics of general form presented in (2.8). This is because, in this case, $\xi = \partial/\partial x^0$ is
not a unit vector field and thus it should be normalized. But this process leads to complicated formulas for
kinematic quantities and Ricci tensor field, which of course makes difficult their study. The question is: Are
there important cosmological models whose metrics have the general form (2.8). The answer is in the affirmative
and it is based on the following two examples. First, the study of the cosmological perturbations of the
FLRW universes is developed with respect to the metric (cf. (10.12) of \cite{ej})

$$ \begin{array}{r} ds^2 = a^2\left\{-(1 + 2\phi)(dx^0)^2 + 2(B_{\vert_i} - S_i)dx^0dx^i\right. \vspace{3mm}\\
\left. + \left[(1 - 2\psi)\gamma_{ij} + 2E_{i{\vert j}}  + 2F_{i{\vert j}}
+ h_{ij}\right]dx^idx^j\right\}.\end{array}\eqno(1.1)$$
Then, we recall that the metric of a Kerr-Newman black hole is given by (cf.(12.3.1 in \cite{w})

$$\begin{array}{r} ds^2 = -\left(\frac{\Delta - a^2(\sin x^2)^2}{\Sigma}\right)(dx^0)^2 - \frac{2a((x^1)^2 + a^2
- \Delta)(\sin x^2)^2}{\Sigma}dx^0dx^3 \vspace{3mm}\\  + \left[\frac{((x^1)^2 + a^2)^2 -
\Delta a^2(\sin x^2)^2}{\Sigma}\right](\sin x^2)^2(dx^3)^2 + \frac{\Sigma}{\Delta}(dx^1)^2
+ \Sigma(dx^2)^2,\end{array}\eqno(1.2)$$
where we put

$$\Delta = (x^1)^2 + a^2 + e^2 - 2mx^1, \ \ \ \Sigma = (x^1)^2 + a^2(\cos x^2)^2. $$
The metric (1.1) was intensively studied with respect to the (1+3) threading of almost FLRW universes. From Chapters
10 and 11 of \cite{ej} we can see that the study is not an easy one in literature. Also, as far we know, very
little has been done with respect to the (1+3) threading theory for the metric (1.2) (cf.[2, 3]).\par
In this paper we present a new approach on the (1+3) threading of spacetime with respect to a congruence of curves
defined by an arbitrary timelike vector field $\xi$. We develop a method that is based on the follwing concepts:\par
(i) \ \ Threading frame and coframe fields.\par
(ii) \ Spatial tensor fields.\par
(iii) Riemannian spatial connection.\newline
The threading frame and coframe fields are naturally constructed from the coordinate fields (cf. (2.3) and (2.4)),
and have a great role throughout the paper. The spatial tensor fields have been used in earlier literature, but
in here we work only with their $3D$ local components with respect to the above special frames (cf.(3.1)). This
brings a substantial simplification into the study of such general metrics. Finally, the Riemannian spatial
connection $\nabla^{\star}$ (cf.(3.12)) is a metric linear connection on the spatial distribution, which
introduces both the spatial and time covariant derivatives. It is important to note that throughout the paper,
all geometric objects and equations involved into the study, are expressed in terms of spatial tensor fields and
their spatial or time covariant derivatives. As the metrics (1.1) and (1.2) fall into the class of the general
metrics given by (2.8), the (1+3) threading theory developed here can be easily applied to their study.\par
Now, we outline the content of the paper. In Section 2 we consider the  orthogonal decomposition (2.1) of the
tangent bundle of the spacetime $(M, g)$, and construct the threading frame and coframe fields
$\{ \partial/\partial x^0, \delta/\delta x^i\}$ and $\{\delta x^0, dx^i\}$, respectively (cf. (2.3) and (2.4)). Also,
we consider the
Riemannian metric $h$ on the spatial distribution $SM$ given by its $3D$ local components (cf. (2.15)). Then,
in Section 3 we introduce the notion of spatial tensor field via $3D$ local components (cf. (3.1)),  and show
that the vorticity, expansion and shear tensor fields $\omega_{ij}, \ \theta_{ij}$ and $\sigma_{ij}$ given by
(3.5a) and (3.8), are indeed spatial tensor fields. Also, we define the Riemannian spatial connection
$\nabla^{\star}$ on $SM$ (cf.(3.11)) and express the Levi-Civita connection $\nabla$ in terms of the local
 coefficients of $\nabla^{\star}$ and the above kinematic quantities (cf. (3.17)). A comparison between the
concepts defined in this paper and the corresponding ones from earlier literature is done in Section 4. In
particular, for a unit timelike vector field we obtain (4.14) for kinematic quantities, and deduce that they
do not depend on the Levi-Civita connection of the spacetime. In Section 5 we express both the curvature tensor
field and the Ricci tensor field  of $(M,g)$ by spatial tensor fields and their spatial and time covariant
derivatives (cf. (5.3), (5.5a), (5.11), (5.12a)). Next, in Section 6 we obtain the Raychaudhuri's equation
(6.1) with respect to an arbitrary timelike vector field, which for a congruence of timelike geodesics takes the
forms (6.5) or (6.28). It is important to note that (6.28) is the main ingredient used in the proof of Lemma
6.2, which should be considered as a completion of Lemma 6.1 that has been the key in the proof of
Penrose-Hawking singularity theorems. Also, we express the non-zero local components of the electric Weyl
curvature tensor field in terms of spatial tensor fields (cf.(6.20)), and deduce new formulas for time
covariant derivatives of the kinematic quantities (cf.(6.9), (6.12), (6.24), (6.25)). Finally, the last three
sections are devoted to the study of a Kerr-Newman black hole via the new approach on the $(1+3)$ threading of
 spacetime developed in the paper. In particular, we characterize spatial geodesics and obtain the $3D$ force
identity (cf.(9.10)). %\newpage

\section{Threading Frame and Coframe Fields}

Let $(M, g)$ be a $4D$ spacetime, and $\xi$ be a timelike vector field  that is globally defined on $M$. Note that
$\xi$ is not necessarily a unit timelike vector field, as it was considered in early literature. The timelike
 congruence determined by $\xi$ is tangent to the fibres of the line bundle $VM$, that we call the {\it time
distribution}. Also, we consider the {\it spatial distribution} $SM$, which is complementary orthogonal to $VM$
in $TM$, that is, we have

$$ TM = VM\oplus SM.\eqno(2.1)$$

Throughout the paper we use the ranges of indices: $i, j, k, ... \in\{1,2,3\}$ and $a, b, c, ...
\in\{0, 1, 2, 3\}$. Also, for any vector bundle $E$ over $M$ denote by $\Gamma(E)$ the ${\cal{F}}(M)$-module of
smooth sections of $E$, where ${\cal{F}}(M)$ is the algebra of smooth functions on $M$.\par
The foliation by curves that is tangent to $VM$, induces a special coordinate system $(x^a)$ such that $\xi =
 {\partial}/{\partial x^0}$. If $(\tilde{x}^a)$ is another coordinate system, then we have

$$\widetilde{x}^i = \widetilde{x}^i(x^1, x^2, x^3); \ \ \ \ \ \ \widetilde{x}^0 = {x}^0 +
f(x^1, x^2, x^3),\eqno(2.2)$$
since $\partial /\partial x^0$ and $\partial /\partial\tilde{x}^0$ represent the same vector field $\xi$, and
hence $\partial\tilde{x}^0 /\partial x^0 = 1.$ Then, from (2.1) we deduce that for each $ \partial /\partial x^i$
 there exist a unique $\delta/\delta x^i \in \Gamma(SM)$ and a unique function $A_i$, such that

$$\frac{\delta }{\delta x^i} = \frac{\partial }{\partial x^i} -
A_i\frac{\partial }{\partial x^0}. \eqno(2.3)$$
This enables us to consider the {\it threading frame field} $\{\partial/\partial x^0, \delta/\delta x^i\}$, and the
{\it threading coframe field} $\{\delta x^0, dx^i\}$, where we put

$$\delta x^0 = dx^0 + A_idx^i.\eqno(2.4)$$
Now, by direct calculations using (2.1)-(2.3), we obtain

$$(a) \ \ \ \frac{\delta}{\delta x^i} = \frac{\partial \widetilde{x}^k}{\partial x^i}
\frac{\delta }{\delta\tilde{x}^k}, \ \ \ (b) \ \ \ \delta\tilde{x}^0 = \delta x^0, \ \ \ (c) \ \ \
A_i = \widetilde{A}_k \frac{\partial \widetilde{x}^k}{\partial x^i} +
\frac{\partial f}{\partial x^i}. \eqno(2.5)$$
Note that $\{\delta/\delta x^i\}$ are transformed exactly as $\{\partial/\partial x^i\}$  on a $3D$ manifold, while
$\{A_i\}$, in general, do not satisfy some $3D$ tensorial transformations.\par
Next, we consider the 1-form $\xi^{\star}$ given by

$$\xi^{\star}(X) = g(X,\xi), \ \ \ \forall X \in \Gamma(TM).\eqno(2.6)$$
The local components of $\xi^{\star}$ with respect to the natural frame field $\{\partial/\partial x^i\}$ are given by

$$\begin{array}{lc}(a) \ \ \ \xi_i = g\left(\frac{\partial }{\partial x^i}, \frac{\partial }{\partial x^0}\right),
\ \ \  (b) \ \ \ \ \xi_0 = g\left(\frac{\partial }{\partial x^0}, \frac{\partial }{\partial x^0} \right) =
-\Phi^2, \end{array}\eqno(2.7)$$
where $\Phi$ is a non-zero function on $M$ which is independent of $x^0$. The above condition on $\Phi$ is
not restrictive for our theory, because most of the important cosmological models satisfy it.\par
According to (2.7), the line element of $g$ is expressed as follows:

$$ds^2 = - \Phi^2(dx^0)^2 + 2\xi_idx^idx^0 + g_{ij}dx^idx^j,\eqno(2.8)$$
where we put

$$g_{ij} = g\left(\frac{\partial }{\partial x^i}, \frac{\partial }{\partial x^j}\right).\eqno(2.9)$$
Taking into account that

$$g(\frac{\delta}{\delta x^i}, \frac{\partial }{\partial x^0}) = 0, \eqno(2.10)$$
and using (2.3) and (2.7), we obtain

$$A_i = -\Phi^{-2}\xi_i,\eqno(2.11)$$
and therefore

$$\frac{\delta }{\delta x^i} = \frac{\partial }{\partial x^i} + \Phi^{-2}\xi_i\frac{\partial }{\partial x^0}.
\eqno(2.12)$$
Multiply (2.5c) by $\Phi^2$, and using (2.11), we obtain

$$\xi_i = \tilde{\xi}_k\frac{\partial\tilde{x}^k}{\partial x^i} - \Phi^2
\frac{\partial f}{\partial x^i}.\eqno(2.13)$$
Hence $\xi_i$ {\it do not define a $3D$ 1-form on $M$}.\par
Now, denote by $h$ the Riemannian metric induced by $g$ on $SM$, and put

$$h_{ij} = h\left(\frac{\delta}{\delta x^i}, \frac{\delta}{\delta x^j}\right) =
g\left(\frac{\delta}{\delta x^i}, \frac{\delta}{\delta x^j}\right).\eqno(2.14)$$
Then by using (2.14), (2.3), (2.11), (2.9) and (2.7), we infer that

$$h_{ij} = g_{ij} + \Phi^{-2}\xi_i\xi_j.\eqno(2.15)$$
Thus $ds^2$ from (2.8) is expressed in terms of the threading coframe field $\{\delta x^0, dx^i\}$ as follows,

$$ds^2 = - \Phi^2(\delta x^0)^2  + h_{ij}dx^idx^j .\eqno(2.16) $$
Note that $h_{ij}$ and the entries $h^{ij}$ of the inverse of the matrix $[h_{ij}]$ are transformed exactly like $3D$
tensor fields, that is, we have

$$(a)  \ \ \ \ h_{ij} = \widetilde{h}_{kh}\frac{\partial \widetilde{x}^{k}}{\partial x^i}
\frac{\partial \widetilde{x}^h}{\partial x^j}, \ \ \ \ (b)  \ \ \ \
 \widetilde{h}^{kh} = h^{ij}\frac{\partial \widetilde{x}^k}{\partial x^i}\frac{\partial \widetilde{x}^h}
{\partial x^j}. \eqno(2.17)$$

\section{ Kinematic Quantities as Spatial Tensor Fields on $(M, g)$}
The purpose of this section is to define the vorticity tensor field, expansion tensor field, expansion scalar
and shear tensor field, as spatial tensor fields on the spacetime $(M, g)$.
First, we give the following definition.\par
A {\it spatial tensor field} $T$ of type $(p, q)$ on $M$, is locally given by $3^{p+q}$ locally defined functions
$T_{i_1 \cdots i_q}^{j_1\cdots j_p}(x)$, satisfying

$$T_{i_1 \cdots i_q}^{j_1\cdots j_p}\frac{\partial\tilde{x}^{k_1}}{\partial x^{j_1}}
\cdots \frac{\partial\tilde{x}^{k_p}}{\partial x^{j_p}}
= \tilde{T}^{k_1 \cdots k_p}_{h_1\cdots h_q}\frac{\partial\tilde{x}^{h_1}}{\partial x^{i_1}}
\cdots \frac{\partial\tilde{x}^{h_q}}{\partial x^{i_q}},\eqno(3.1)$$
with respect to the transformations (2.2). In other words, the local components of a spatial tensor field on $M$
should satisfy the same transformations as the local components of a tensor field on a 3-dimensional manifold.
From (2.17) we see that $h_{ij}$ (resp. $h^{ij}$) define a spatial tensor field of type (0,2) (resp.(2,0)) on $M$.\par
By using (2.5a) and taking into account that

$$\frac{\partial\Phi}{\partial x^0} = 0, \eqno(3.2)$$
we deduce that

$$c_i = \Phi^{-1}\frac{\delta\Phi}{\delta x^i} = \Phi^{-1}\frac{\partial\Phi}{\partial x^i}, \eqno(3.3)$$
define a spatial tensor field of type (0,1). Next, by direct calculations using (2.3), (2.11) and (3.3), we deduce that

$$(a) \ \ \left [\frac{\delta }{\delta x^j}, \frac{\delta }{\delta x^i}\right] =
2\omega_{ij}\frac{\partial }{\partial x^0},
\ \ (b) \ \ \left [\frac{\partial }{\partial x^0}, \frac{\delta }{\delta x^i}\right] =
a_i\frac{\partial }{\partial x^0},\eqno(3.4)$$
where we put

$$\begin{array}{lc}
(a) \ \ \ \omega_{ij} = \frac{1}{2}\left\{\frac{\delta A_j}{\delta x^i} - \frac{\delta A_i}{\delta x^j}\right\}
= \Phi^{-2}\left\{c_i\xi_j - c_j\xi_i + \frac{1}{2}\left(\frac{\delta\xi_i}{\delta x^j} - \frac{\delta\xi_j}{\delta x^i}
\right)\right\}, \vspace{2mm} \\ (b) \ \ \
 a_i = -\frac{\partial A_i}{\partial x^0} = \Phi^{-2}\frac{\partial\xi_i}{\partial x^0}. \end{array}\eqno(3.5)$$
 Now, apply $\delta/\delta x^j$ and $\partial/\partial x^0$ to (2.5c) and by using (2.5a) and (2.5b), we infer that

 $$\begin{array}{c}(a)\ \ \ \frac{\delta A_i}{\delta x^j} = \frac{\delta\widetilde{A}_k}{\delta\widetilde{x}^h}
 \frac{\partial\widetilde{x}^h}{\partial x^j}\frac{\partial\widetilde{x}^k}{\partial x^i}
+ \widetilde{A}_k\frac{\partial^2\tilde{x}^k}{\partial x^i\partial x^j}
+  \frac{\partial^2 f} {\partial x^i\partial x^j}, \vspace{3mm} \\
(b) \ \ \ \frac{\partial A_i}{\partial x^0} = \frac{\partial \widetilde{A}_k}{\partial\tilde{x}^0}
\frac{\partial\tilde{x}^k}{\partial x^i}.\end{array} \eqno(3.6)$$
 Then, by using (3.6) into (3.5), we obtain

 $$\begin{array}{c}(a)\ \ \ \omega_{ij} = \widetilde{\omega}_{kh}
 \frac{\partial\widetilde{x}^k}{\partial x^i}\frac{\partial\widetilde{x}^h}{\partial x^j}, \ \ \
(b) \ \ \ a_i  =  \widetilde{a}_k\frac{\partial\tilde{x}^k}{\partial x^i}.\end{array} \eqno(3.7)$$
Hence $\omega_{ij}$ and $a_i$ define spatial tensor fields of type (0.2) and (0,1), respectively. We call
$\omega = (\omega_{ij})$ the {\it vorticity tensor field} for the timelike congruence defined by $\xi$ on $M$.
From (3.4a) we see that the spatial distribution $SM$ is integrable if and only if $\omega_{ij}$ vanish
identically on $M$.\par
 Next, we define

$$(a) \ \ \Theta_{ij} = \frac{1}{2}\frac{\partial h_{ij}}{\partial x^0}, \ \ \ (b) \ \ \
 \Theta = h^{ij}\Theta_{ij}, \ \ \ (c) \ \ \ \sigma_{ij} = \Theta_{ij} - \frac{1}{3}\Theta h_{ij}. \eqno(3.8)$$
Then, take derivatives with respect to $x^0$ in (2.17a), and obtain

$$ \Theta_{ij} = \widetilde{\Theta}_{kh}\frac{\partial\widetilde{x}^k}{\partial x^i}
 \frac{\partial\widetilde{x}^h}{\partial x^j}, $$
that is, $ \Theta_{ij}$ define a spatial tensor field of type (0,2). We call it the {\it expansion tensor field}.
 Clearly $\Theta$ from (3.8b) is a function, and $\sigma_{ij}$ from (3.8c) define a trace-free spatial tensor
field of type (0,2). We call $\Theta$ the {\it expansion scalar} and $\sigma_{ij}$ the {\it shear tensor filed} for
the congruence. According to the terminology from earlier literature, we call $\{\omega_{ij}, \Theta_{ij},
\Theta, \sigma_{ij}\}$ given by (3.5a) and (3.8), the {\it kinematic quantities} with respect to the congruence
of curves defined by the timelike vector field $\xi = \partial/\partial x^0$.\par
Raising and lowering latin indices is done by using $h^{ij}$ and $h_{ij}$, as follows

$$(a) \ \ \ \omega_j^k = h^{ki}\omega_{ij}, \ \ \ (b) \ \ \ \omega^{kh} = h^{ki}h^{hj}\omega_{ij}. \eqno(3.9)$$

In order to define covariant derivatives of the above kinematic quantities, we consider the Levi-Civita
connection $\nabla$ on $(M, g)$ given by (cf. \cite{o}, p.61)

$$\begin{array}{lc}
2g(\nabla_XY, Z) =  X(g(Y, Z)) + Y(g(Z, X)) - Z(g(X, Y)) \vspace{4mm}
\\+ g([X, Y], Z) -  g([Y, Z], X) +   g([Z, X], Y),\end{array}\eqno(3.10) $$
for all $X, Y, Z \in \Gamma(TM)$. Then define the linear connection $\nabla^{\star}$ on the spatial distribution
as the spatial projection of $\nabla$ on $SM$, that is, we have

$$(a) \ \ \ \nabla^{\star}_{X}sY = s{\nabla}_XsY, \ \ \ \forall \ \ X,Y \in(\Gamma(SM),\eqno(3.11)$$
where $s$ is the projection morphism of $TM$ on $SM$ with respect to (2.1). Note that $\nabla^{\star}$ is a
metric linear connection on $SM$. We call it the {\it Riemannian spatial connection}.\par
{\bf Remark 3.1} The Riemannian spatial connection $\nabla^{\star}$ is different from the three-dimensional operator
$\bar{\nabla}$ that has been used in earlier literature (cf. (4.19) of [1]). $\nabla^{\star}$ is a linear connection
on $SM$ and therefore defines covariant derivatives of any spatial tensor field with respect to vector fields on
$M$. On the contrary, $\bar{\nabla}$ is an operator which acts on tensor fields on $M$, but in general, does not define
a linear connection on $M$.\ $\Box$ \par
 Locally, we put

$$(a) \ \ \nabla^{\star}_{\frac{\delta }{\delta x^j}}\frac{\delta }{\delta x^i} =
\Gamma^{\star \;k}_{i\ \;j}\frac{\delta }{\delta x^k}, \ \ \ (b) \ \ \ \
\nabla^{\star}_{\frac{\partial }{\partial x^0}}\frac{\delta }{\delta x^i} =
\Gamma^{\star \;k}_{i \ \;0}\frac{\delta }{\delta x^k}.\eqno(3.12) $$
Then, take $X = \delta/\delta x^j, Y = \delta/\delta x^i$ and $Z = \delta/\delta x^h$ in (3.10) and by using
(3.11), (3.12a), (2.14) and (3.4a), we obtain

$$\Gamma^{\star \;k}_{i\ \;j} = \frac{1}{2}h^{kh}\left\{\frac{\delta h_{hj}}{\delta x^i} +
\frac{\delta h_{hi}}{\delta x^j} - \frac{\delta h_{ij}}{\delta x^h}\right\}. \eqno(3.13) $$
Similarly, we deduce that

$$\Gamma^{\star \;k}_{i \ \;0} = \Theta_i^k + \Phi^2\omega_i^k.\eqno(3.14) $$

Now, consider a spatial tensor field $T$ of type $(p,q)$. Then  $\nabla^{\star}_{\frac{\delta}{\delta x^k}}T$
and $\nabla^{\star}_{\frac{\partial}{\partial x^0}}T$ are spatial tensor of type $(p, q+1)$ and
$(p, q)$, respectively. As an example, we consider  $T = (T^i_j)$, and obtain

$$ \begin{array}{c}(a) \ \ \
T_{j\vert_{k}}^{i} = \frac{\delta T_j^i}{\delta x^k} + T_j^h\Gamma^{\star \;i}_{h\ \;k}
  - T_{h}^{i}\Gamma^{\star \;h}_{j\ \;k},\vspace{3mm}\\ (b) \ \ \
 T_{j\vert_{0}}^i = \frac{\partial T_j^i}{\partial x^0} + T_j^h\Gamma^{\star \;i}_{h\ \;0}
  - T_h^i\Gamma^{\star \;h}_{j\ \;0}.\end{array}  \eqno(3.15)$$
We call (3.15a) (resp. (3.15b)) the {\it spatial} (resp. {\it time}) {\it covariant derivative} of $T$.
As $\nabla^{\star}$ is a metric connection on $SM$, we have

$$(a) \ \ \ h_{ij\vert_{k}} = 0, \ \ \ (b) \ \ \ h^{ij}_{\ \ \vert_{k}}
= 0,\ \ \ (c) \ \ \ h_{ij\vert_{0}} = 0, \ \ \ (d) \ \ \ h^{ij}_{\ \ \vert_{0}}
= 0.\eqno(3.16)$$

Finally, by using (3.10), the above spatial tensor fields and the local coefficients of $\nabla^{\star}$,
we express  the Levi-Civita connection $\nabla$ on $(M, g)$ as follows

$$\begin{array}{lc}
(a) \ \ \nabla_{\frac{\delta }{\delta x^j}}\frac{\delta }{\delta x^i} =  \Gamma^{\star \;k}_{i\ \;j}
\frac{\delta }{\delta x^k} + \left(\omega_{ij} + \Phi^{-2}\Theta_{ij}\right)\frac{\partial }{\partial x^0},
\vspace{4mm}\\
(b) \ \ \nabla_{\frac{\partial }{\partial x^0}}\frac{\delta }{\delta x^i} = \left(\Theta_i^k
+ \Phi^2\omega_i^k\right)\frac{\delta }{\delta x^k} + \left(a_i + c_i\right)
\frac{\partial }{\partial x^0}, \vspace{4mm}\\
(c) \ \ \nabla_{\frac{\delta }{\delta x^i}}\frac{\partial }{\partial x^0} =
\left(\Theta_i^k + \Phi^2\omega_i^k\right)\frac{\delta }{\delta x^k} + c_i\frac{\partial }{\partial x^0},
\vspace{4mm}\\
(d) \ \ \nabla_{\frac{\partial }{\partial x^0}}\frac{\partial }{\partial x^0} =
 \Phi^2\left(a^k + c^k\right)\frac{\delta }{\delta x^k} . \end{array}\eqno(3.17)$$

\section{Comparison with Concepts from Earlier Literature}

In the previous section we introduced the kinematic quantities on a spacetime $(M, g)$ with respect to the
congruence that is tangent to an arbitrary timelike vector field $\xi$. If in particular, $\xi$ is a unit
timelike vector field, the configuration of the spacetime with respect to the congruence of timelike curves
determined by $\xi$ is known in literature as (1+3) {\it threading of spacetime} (cf.[1,5]). In
this section we show that for $\Phi^2 = 1$ in (3.5) and (3.8) (that is, $\xi$ is a unit vector field), we obtain
the well known kinematic quantities from earlier literature.\par
First, by using (2.6) and taking into account that $\nabla$ is a metric connection, we obtain

$$(\nabla_X\xi^{\star})(Y) = g(Y, \nabla_X\xi), \ \ \ \forall \ X, Y \in \Gamma(TM). \eqno(4.1)$$
 Then, consider the threading frame $\{\partial/\partial x^0, \delta/\delta x^i\}$ and using (3.17), we infer that

$$\begin{array}{l}(a) \ \ \
\left(\nabla_{\frac{\delta}{\delta x^j}}\xi^{\star}\right)\left(\frac{\delta}{\delta x^i}\right)
= \Theta_{ij} + \Phi^2\omega_{ij}, \vspace{3mm} \\ (b) \ \ \
\left(\nabla_{\frac{\partial }{\partial x^0}}\xi^{\star}\right)\left(\frac{\delta}{\delta x^i}\right)
= \Phi^2(a_i + c_i), \vspace{3mm} \\
(c) \ \ \
\left(\nabla_{\frac{\delta}{\delta x^i}}\xi^{\star}\right)\left(\frac{\partial}{\partial x^0}\right)
= -\Phi^2c_i, \vspace{3mm} \\
(d) \ \ \
\left(\nabla_{\frac{\partial}{\partial x^0}}\xi^{\star}\right)\left(\frac{\partial}{\partial x^0}\right)
= 0. \end{array}\eqno(4.2)$$
Next, taking into account (2.12) and (2.7b), we express the natural frame field as follows

$$\frac{\partial }{\partial x^a} = \delta^i_a\frac{\delta }{\delta x^i} - \Phi^{-2}\xi_a
\frac{\partial }{\partial x^0}, \ \ \ a \in\{0, 1, 2, 3\}.\eqno(4.3)$$
Then, consider the {\it covariant acceleration vector field}

$${\stackrel{.}{\xi}}_a = \left(\nabla_{\frac{\partial}{\partial x^0}}\xi^{\star}\right)
 \left(\frac{\partial}{\partial x^a}\right), \eqno(4.4)$$
and using (4.2b) and (4.2d), we obtain

$${\stackrel{.}{\xi}}_a = \Phi^2\delta_a^ib_i, \eqno(4.5)$$
where we put

$$b_i = a_i + c_i.\eqno(4.6) $$
Thus, we deduce that the congruence defined by $\xi$ {\it is a congruence of timelike geodesics, if and
only if, we have}

$$b_i = 0, \forall \ i \in \{1, 2, 3\}.\eqno(4.7)$$
For this reason we call $b_i$ the {\it geodesic spatial covector field} of the congruence.\par
Now, by direct calculations, using (4.3) and (4.2), we infer that

$$\begin{array}{l}(a) \ \ \
\left(\nabla_{\frac{\partial}{\partial x^b}}\xi^{\star}\right)\left(\frac{\delta}{\delta x^i}\right)
= \delta_b^j\left(\Theta_{ij} + \Phi^2\omega_{ij}\right) - b_i\xi_b, \vspace{3mm} \\ (b) \ \ \
\left(\nabla_{\frac{\partial}{\partial x^b}}\xi^{\star}\right)\left(\frac{\partial}{\partial x^0}\right)
= -\Phi^2\delta_b^jc_j.\end{array}\eqno(4.8)$$
Taking into account (4.3), (4.8) and (4.5), we find

$$\nabla_b\xi_a = - \Phi^{-2}\xi_b{\stackrel{.}{\xi}}_a + \delta_a^i\delta_b^j\left(\Theta_{ij} +
\Phi^2\omega_{ij}\right) + \xi_a\delta_b^jc_j.\eqno(4.9)$$
Now, we suppose that $\xi$ is a unit vector field. Then, according to the formula (4.38) from \cite{ej}, p.85, we have

$$\nabla_b\xi_a = - \xi_b{\stackrel{.}{\xi}}_a + \sigma_{ab} + \frac{1}{3}\Theta h_{ab} + \omega_{ab},
\eqno(4.10)$$
where $\sigma_{ab}, \theta, h_{ab}$ and $\omega_{ab}$ are quantities defined in earlier literature. On the other
hand, in this case we have $\Phi^2 = 1,$ and from (4.9) we obtain

$$\nabla_b\xi_a = -\xi_b{\stackrel{.}{\xi}}_a + \delta_a^i\delta_b^j\left(\Theta_{ij} +
\omega_{ij}\right),\eqno(4.11)$$
since $c_j = 0,$ for all $j \in \{1, 2, 3\}$.\par
Comparing the symmetric and skew-symmetric parts in (4.10) and (4.11) we deduce that

$$\begin{array}{l} (a) \ \ \ \sigma_{ab} = \delta_a^i\delta_b^j\Theta_{ij} - \frac{1}{3}\Theta h_{ab}, \vspace{3mm} \\
 (b) \ \ \ \omega_{ab} = \delta_a^i\delta_b^j\omega_{ij}.\end{array}\eqno(4.12)$$
Finally, comparing (4.12a) with (4.31) from \cite{ej}, p.81, we obtain

$$\Theta_{ab} = \delta_a^i\delta_b^j\Theta_{ij}.\eqno(4.13)$$
According to (3.8), (3.5a), (4.12) and (4.13), we conclude that in case $\xi$ is a unit vector field, the
only possible non-zero local components of expansion, shear and vorticity tensor fields from earlier literature
are the following

$$\begin{array}{l}(a) \ \ \
\Theta_{ij} =  \frac{1}{2}\frac{\partial h_{ij}}{\partial x^0}, \ \ \ (b) \ \ \
\sigma_{ij} = \frac{1}{2}\frac{\partial h_{ij}}{\partial x^0} - \frac{1}{3}\Theta h_{ij},\vspace{3mm} \\
(c) \ \ \ \omega_{ij} = \frac{1}{2}\left\{\frac{\delta\xi_i}{\delta x^j} -
\frac{\delta\xi_j}{\delta x^i}\right\}.\end{array} \eqno(4.14)$$
As far as we know, the formulas (4.14) do not appear in earlier literature. Due to them we can state that {\it
the expansion, shear and vorticity tensor fields do not depend on the Levi-Civita connection of the spacetime}
$(M, g)$. Of course, due to (3.5a) and (3.8), this conclusion is still valid for the general case of a
congruence defined by an arbitrary timelike vector field $\xi$.\par
We close this section with an interesting property of vorticity tensor field. Suppose that we have a congruence
of geodesics defined by $\xi$. Then, according to (4.7), (4.6), (3.3) and (3.5b), we have

$$\frac{\partial \xi_i}{\partial x^0} + \Phi\frac{\partial \Phi}{\partial x^i} = 0.\eqno(4.15)$$
By using (4.15), (3.2) and (3.3), we infer that

$$\begin{array}{c}(a) \ \ \ \frac{\partial^2\xi_i}{(\partial x^0)^2} = 0, \ \ \ (b) \ \
\ c_i\frac{\partial\xi_j}{\partial x^0} = c_j\frac{\partial\xi_i}{\partial x^0}, \ \ \ (c) \ \ \
\frac{\partial^2\xi_i}{\partial x^j\partial x^0} = \frac{\partial^2\xi_j}{\partial x^i\partial x^0}.
\end{array}\eqno(4.16)$$
Then, take the derivative with respect to $x^0$ in (3.5a) and by (3.2), (3.3) and (4.16), we obtain

$$\frac{\partial \omega_{ij}}{\partial x^0} = 0.\eqno(4.17)$$
Thus, we can state that {\it the vorticity tensor field for a congruence of timelike geodesics of a spacetime,
is independent of time}. In particular, this is true for a congruence of geodesics with respect to a unit
vector field $\xi$. However, we did not see this result in earlier literature. This is because the formulas
(3.5a), (3.8) and (4.14) we deduced for the kinematic quantities are much simpler than the ones by means of
Levi-Civita connection.

\section{Curvature and Ricci Tensor Fields of a Spacetime via Spatial Tensor Fields}

In this section we show that the curvature tensor field of $(M, g)$ is completely determined by three spatial
tensor fields $R_{ijkh}, R_{i0kh}$ and $R_{i0k0}$ (cf. (5.3), (5.5a)). A similar result we obtain for the Ricci
tensor of $(M, g)$ (cf. (5.11), (5.12a)). Note that all these spatial tensor fields are expressed in terms of
the curvature and Ricci tensor fields of the of the Riemannian spatial connection, and of all kinematic
quantities introduced in Section 3.\par
In what it follows, $R$ denotes both the curvature tensor field of $(M, g)$ of type (0,4) and (1,3), given by

$$\begin{array}{l} \vspace{2mm}\\
(a) \ \ \ R(X, Y, Z, U) = g(R(X, Y, U), Z),\vspace{2mm}\\
(b) \ \ \ R(X, Y, U) = \nabla_X\nabla_YU - \nabla_Y\nabla_XU - \nabla_{[X, Y]}U,\end{array}\eqno(5.1)$$
for all $X, Y, Z, U \in \Gamma(TM)$. Then the curvature tensor field of $(M, g)$ is completely determined by its
local components

$$\begin{array}{l}(a) \ \ \ R_{ijkh} = R(\frac{\delta }{\delta x^h}, \frac{\delta }{\delta
 x^k}, \frac{\delta }{\delta x^j}, \frac{\delta }{\delta x^i}), \vspace{3mm}\\
 (b) \ \ \ R_{i 0 kh} = R(\frac{\delta }{\delta x^h}, \frac{\delta }{\delta x^k},
\frac{\partial }{\partial x^0}, \frac{\delta }{\delta x^i}), \vspace{3mm}\\
(c) \ \ \ R_{i 0 k 0} = R(\frac{\partial }{\partial x^0}, \frac{\delta }{\delta x^k},
\frac{\partial }{\partial x^0}, \frac{\delta }{\delta x^i}).\end{array}\eqno(5.2)$$
By direct calculations, using (5.2), (5.1), (3.17), (3.4), (3.3) and (4.6), we obtain

$$\begin{array}{lc}(a) \ \  \ R_{ijkh} = R^{\star}_{ijkh} + \omega_{ik}\Theta_{jh} - \omega_{ih}\Theta_{jk} +
\Phi^{-2}\left(\Theta_{ik}\Theta_{jh} - \Theta_{ih}\Theta_{jk}\right)\vspace{3mm}\\
\hspace*{22mm}+ \Phi^2\left(\omega_{ik}\omega_{jh}
- \omega_{ih}\omega_{jk}\right) + \Theta_{ik}\omega_{jh} - \Theta_{ih}\omega_{jk}, \vspace{3mm}\\
(b) \ \ \ R_{i0kh} = \Theta_{ih\vert_k} -  \Theta_{ik\vert_h} + \Theta_{ik}c_h - \Theta_{ih}c_k
\vspace{3mm}\\ \hspace*{22mm}+ \Phi^2\left\{\omega_{ih\vert_k} -  \omega_{ik\vert_h} + \omega_{ih}c_k -
\omega_{ik}c_h + 2\omega_{kh}b_i \right\}, \vspace{3mm}\\
(c) \ \ \ R_{i0k0} = \Phi^2\left\{b_ib_k + b_{i\vert_k} + \omega_{kh}\Theta_i^h
- \omega_{ih}\Theta_k^h - \omega_{ik\vert_0} - \Phi^2\omega_{ih}\omega_k^h\right\}\vspace{3mm}\\
\hspace*{22mm}- \Theta_{ik\vert_0} - \Theta_{ih}\Theta_k^h, \end{array}\eqno(5.3)$$
where $R^{\star}_{ijkh}$ are the local components of the curvature tensor field of the Riemannian spatial
connection defined as in (5.2a), and given by

$$\begin{array}{lc}R^{\star}_{ijkh} = h_{jl}\left\{\frac{\delta\Gamma_{i\;\;k}^{\star l}}{\delta x^h} -
\frac{\delta\Gamma_{i\;\;h}^{\star l}}{\delta x^k} + \Gamma_{i\;\;k}^{\star n}\Gamma_{n\;\;h}^{\star l}
- \Gamma_{i\;\;h}^{\star n}\Gamma_{n\;\;k}^{\star l}  \right.\vspace{3mm}\\ \left. \hspace*{12mm}
- 2\omega_{kh}\left(\Theta_i^l + \Phi^2\omega_i^l\right)\right\}.\end{array}\eqno(5.4)$$
Taking the symmetric and skew-symmetric parts, in (5.3c) we deduce that

$$\begin{array}{l} (a) \ \ \ R_{i0k0} = \Phi^2\left\{b_ib_k + \frac{1}{2}\left(b_{i\vert_k} + b_{k\vert_i}\right)
 - \Phi^2\omega_{ih}\omega_k^h\right\}\vspace{3mm}\\
\hspace*{22mm}- \Theta_{ik\vert_0} - \Theta_{ih}\Theta_k^h,\vspace{3mm}\\
(b) \ \ \ \omega_{ik\vert_0} = \omega_{kh}\Theta_i^h  - \omega_{ih}\Theta_k^h +
\frac{1}{2}\left(b_{i\vert_k} - b_{k\vert_i}\right). \end{array}\eqno(5.5)$$
Now, consider an orthonormal frame field $\{E_k, \Phi^{-1}\frac{\partial }{\partial x^0}\}$ on $M$ and put

$$E_k = E_k^i\frac{\delta }{\delta x^i}.\eqno(5.6)$$
Then, we deduce that

$$h^{ij} =  \sum_{k =1}^3E_k^iE_k^j.\eqno(5.7)$$
According to  \cite{o}, p.87, the Ricci tensor of $(M, g)$ is given by

$$Ric(X, Y) = \sum_{k =1}^3R(E_k, X, E_k, Y) - \Phi^{-2}R(\frac{\partial }{\partial x^0}, X,
\frac{\partial }{\partial x^0}, Y). \eqno(5.8)$$
Then, by using (5.8), (5.6), (5.7) and (5.2), we obtain

$$\begin{array}{c}(a) \ \
 R_{ik} =  h^{jh}R_{ijkh} -  \Phi^{-2}R_{i 0 k 0}, \ \ \ (b) \ \ R_{i 0} = h^{jh}R_{j 0 hi},
\vspace{4mm}\\ \ \ (c) \ \ \ R_{00} =  h^{jh}R_{j 0 h 0},\end{array}\eqno(5.9) $$
where we put

$$\begin{array}{lc}(a) \ \ \  R_{ik} =
 Ric\left(\frac{\delta }{\delta x^k}, \frac{\delta }{\delta x^i}\right), \ \ \ \ (b) \ \ \ R_{i0} =
 Ric\left(\frac{\partial }{\partial x^0}, \frac{\delta }{\delta  x^i}\right),
\vspace{4mm}\\ \ \ (c) \ \ \  R_{00} =  Ric\left(\frac{\partial }{\partial x^0},
\frac{\partial }{\partial x^0}\right).\end{array}\eqno(5.10) $$
By using (5.3a), (5.3b) and (5.5a) into (5.9), we deduce that

$$\begin{array}{l} (a) \ \ \ R_{ik} = R_{i\ \ kh}^{\star h} + \Phi^{-2}\left(\Theta_{ik\vert_{0}}
 + \Theta\Theta_{ik}\right) - b_ib_k - \frac{1}{2}\left(b_{i\vert_k} + b_{k\vert_i}\right)\vspace{3mm}\\
\hspace*{22mm}+ \Theta\omega_{ik} + \omega_{kh}\Theta_i^h - \omega_{ih}\Theta_k^h,  \vspace{3mm}\\
(b) \ \ \ R_{i0} = \Theta^k_{i\vert_k} - \Theta_{\vert_i} + \Theta c_i - \Theta_i^kc_k \vspace{3mm}\\
\hspace*{22mm} +
\Phi^{2}\left\{\omega^k_{i\ \vert_k} + \omega_i^kc_k + 2\omega_i^kb_k \right\},\vspace{3mm}\\
(c) \ \ \ R_{00} = \Phi^2\left\{b_kb^k + b^k_{\;\vert_k} + \Phi^2\omega_{kh}\omega^{kh}\right\} -
\Theta_{\vert_0} - \Theta_{kh}\Theta^{kh}, \end{array} \eqno(5.11)$$
 where $\Theta^k_{i\vert_k}$, $\omega^k_{i\ \vert_k}$ and $b^k_{\ \vert_k}$ are {\it spatial divergences} given
 by formulas deduced from (3.15a). Now, take symmetric and skew-symmetric parts in (5.11a) and obtain

$$\begin{array}{l}(a) \ \ \
R_{ik} = R^{\star}_{ik} + \Phi^{-2}\left(\Theta_{ik\vert_{0}} + \Theta\Theta_{ik}\right) -
b_ib_k - \frac{1}{2}\left(b_{i\vert_{k}} + b_{k\vert_{i}}\right),\vspace{2mm}\\
(b)\ \ \ \frac{1}{2}\left(R_{i\ kh}^{\star h} - R_{k\ ih}^{\star h}\right) = \omega_{ih}\Theta_k^h -
\omega_{kh}\Theta_i^h - \Theta\omega_{ik},
\end{array}\eqno(5.12)$$
 where we put

 $$R^{\star}_{ik} = \frac{1}{2}\left(R_{i\ kh}^{\star h} + R_{k\ ih}^{\star h}\right). \eqno(5.13)$$
We call $R^{\star}_{ik}$ the {\it spatial Ricci tensor} of the  spacetime $(M, g)$.\par
 From (5.12b) we see that if the spatial distribution is integrable,  then we have

$$R_{i\ kh}^{\star h} = R_{k\ ih}^{\star h}.\eqno(5.14) $$
In this case, we have

$$R^{\star}_{ik} = R_{i\ kh}^{\star h}. \eqno(5.15)$$
Also, note that because in this particular case the vorticity vanishes identically, from (5.3a), (5.3b),
(5.5a), (5.11a), (5.11b) and (5.12a), we deduce that the curvature and Ricci tensors on $(M, g)$ are
expressed as follows

$$\begin{array}{lc}(a) \ \  \ R_{ijkh} = R^{\star}_{ijkh} +
\Phi^{-2}\left(\Theta_{ik}\Theta_{jh} - \Theta_{ih}\Theta_{jk}\right),\vspace{3mm}\\
(b) \ \ \ R_{i0kh} = \Theta_{ih\vert_k} - \Theta_{ik\vert_h} + \Theta_{ik}c_h - \Theta_{ih}c_k,
\vspace{3mm}\\
(c) \ \ \ R_{i0k0} = \Phi^2\left\{b_ib_k + \frac{1}{2}\left(b_{i\vert_k} + b_{k\vert_i}\right)\right\}
- \Theta_{ik\vert_0} - \Theta_{ih}\Theta_k^h, \end{array}\eqno(5.16)$$
and

$$\begin{array}{l} (a) \ \ \ R_{ik} = R_{ik}^{\star} + \Phi^{-2}\left(\Theta_{ik\vert_{0}}
 + \Theta\Theta_{ik}\right) - b_ib_k - \frac{1}{2}\left(b_{i\vert_k} + b_{k\vert_i}\right),\vspace{3mm}\\
(b) \ \ \ R_{i0} = \Theta^k_{i\vert_k} - \Theta_{\vert_i} + \Theta c_i - \Theta_i^kc_k, \vspace{3mm}\\
(c) \ \ \ R_{00} = \Phi^2\left\{b_hb^h + b^h_{\;\vert_h}\right\} -
\Theta_{\vert_0} - \Theta_{kh}\Theta^{kh}, \end{array} \eqno(5.17)$$
where $R_{ik}^{\star}$ is given by (5.15).

\section{Raychaudhuri's Equation and Time Covariant Derivatives of Kinematic Quantities}

First, by using (3.8b) and (3.8c) into (5.11c), we infer

$$\Theta_{\vert_0} = \Phi^4\omega^2 - \sigma^2 - \frac{1}{3}\Theta^2 + \Phi^2\left(b^h_{\ \vert_h} + b^2\right)
 - R_{00},\eqno(6.1) $$
where we put

$$(a) \ \ \ \omega^2 = \omega_{kh}\omega^{kh}, \ \ \ (b) \ \ \ \sigma^2 = \sigma_{kh}\sigma^{kh}, \ \ \ (c)
\ \ \ b^2 = b_hb^h. \eqno(6.2)$$
In particular, if $\xi$ is a unit timelike vector field, that is

$$ \Phi^2 = 1,\eqno(6.3)$$
(6.1) becomes

$$\Theta_{\vert_0} = \omega^2 - \sigma^2 - \frac{1}{3}\Theta^2 + b^h_{\ \vert_h} + b^2
 - R_{00},\eqno(6.4) $$
which is Raychaudhuri's equation expressed in terms of local components of spatial tensor fields introduced in
the present paper. Thus, we are entitled to call (6.1) the {\it generalized Raychaudhuri's equation} with respect
to a congruence defined by an arbitrary timelike vector field $\xi$.\par
According to (4.7), in case of a timelike congruence of geodesics, (6.1) and (6.4) become

$$\Theta_{\vert_0} = \Phi^4\omega^2 - \sigma^2 - \frac{1}{3}\Theta^2  - R_{00},\eqno(6.5) $$
and

$$\Theta_{\vert_0} = \omega^2 - \sigma^2 - \frac{1}{3}\Theta^2  - R_{00},\eqno(6.6) $$
respectively.\newline
{\bf Remark 6.1.} Formally, (6.6) looks like (9.2.11) in \cite{w}, but we should note that $\omega, \sigma$ and
$\Theta$ from (6.6) are calculated via their $3D$ spatial components (see  (3.5a), (3.8) and (6.2)), while in
 \cite{w} they are calculated in terms of the $4D$ local components with respect to the natural frame field
$\{\partial/\partial x^a\}$. \ $\Box$ \par
Next, observe that (5.5b) gives a formula for the time covariant derivative of vorticity tensor field. By using
(3.8b) and (3.8c) into (5.5b) we find

$$\omega_{ik\vert_0} = \omega_{kh}\sigma_i^h  - \omega_{ih}\sigma_k^h - \frac{2}{3}\Theta\omega_{ik}
+ \frac{1}{2}\left(b_{i\vert_k} - b_{k\vert_i}\right). \eqno(6.7)$$
Now, from (5.5a) we deduce that the time covariant derivative of the expansion tensor field is given by

$$\Theta_{ik\vert_0} = \Phi^2\left\{b_ib_k + \frac{1}{2}\left(b_{i\vert_k} + b_{k\vert_i}\right)
 - \Phi^2\omega_{ih}\omega_k^h\right\} - \Theta _{ih}\Theta_k^h - R_{i0k0}.\eqno(6.8)$$
Another formula in terms of Ricci tensors is deduced from (5.12a):

$$\Theta_{ik\vert_0} = -\Theta\Theta_{ik} + \Phi^2\left\{b_ib_k + \frac{1}{2}\left(b_{i\vert_{k}} +
b_{k\vert_{i}}\right) + R_{ik} - R_{ik}^{\star}\right\}.\eqno(6.9)$$
Taking time covariant derivative in (3.8c), and using (3.16c), (6.8) and (6.1), we infer that

$$\begin{array}{l}
\sigma_{ik\vert_0} = \Phi^2\left\{b_ib_k + \frac{1}{2}\left(b_{i\vert_{k}} + b_{k\vert_{i}}\right)
- \frac{1}{3}\left(b^h_{\ \vert_h} + b^2\right)h_{ik}\right.\vspace{2mm}\\ \left.- \Phi^2\left(\omega_{ih}\omega_k^h +
\frac{1}{3}\omega^2h_{ik}\right) \right\} + \frac{1}{3}\sigma^2h_{ik} - \sigma_{ih}\sigma^h_k - \frac{2}{3}
\Theta\sigma_{ik} - \widetilde{R}_{i0k0},\end{array}\eqno(6.10)$$
where $\widetilde{R}_{i0k0}$ is the trace-free part of the spatial tensor field $R_{i0k0}$, given by

$$ \widetilde{R}_{i0k0} = R_{i0k0} -\frac{1}{3}R_{00}h_{ik}.\eqno(6.11)$$
In a similar way, but using (6.9) instead of (6.8), we obtain

$$
\begin{array}{l}\sigma_{ik\vert_0} = -\Theta\sigma_{ik} + \frac{1}{3}\left(\sigma^2 - \frac{2}{3}\Theta^2
+ R_{00}\right)h_{ik}+ \Phi^2\left\{b_ib_k\right.\vspace{2mm}\\ \left. + \frac{1}{2}\left(b_{i\vert_{k}}
  + b_{k\vert_{i}}\right) - \frac{1}{3}\left(b^h_{\ \vert_h} + b^2 + \Phi^2\omega^2\right)h_{ik}+
R_{ik} - R_{ik}^{\star}\right\}.\end{array}\eqno(6.12)$$
Now consider the Weyl tensor field in $(M, g)$ given by

$$\begin{array}{l} C_{abcd} = \bar{R}_{abcd} + \frac{1}{2}\left\{g_{ad}\bar{R}_{bc} + g_{bc}\bar{R}_{ad}
- g_{ac}\bar{R}_{bd} - g_{bd}\bar{R}_{ac}\right\}\vspace{2mm}\\ \hspace*{12mm} +
\frac{1}{6}{\bf R}\left\{g_{ac}g_{bd} - g_{ad}g_{bc}\right\},\end{array}\eqno(6.13)$$
where we put

$$\begin{array}{l}(a) \ \ \ C_{abcd} = C\left(\frac{\partial }{\partial x^d}, \frac{\partial }{\partial
 x^c}, \frac{\partial }{\partial x^b}, \frac{\partial }{\partial x^a}\right), \vspace{3mm}\\
 (b) \ \ \ \bar{R}_{abcd} = R(\frac{\partial }{\partial x^d}, \frac{\partial }{\partial x^c},
\frac{\partial }{\partial x^b}, \frac{\partial }{\partial x^a}) \vspace{3mm}\\
(c) \ \ \ \bar{R}_{ab} = Ric(\frac{\partial }{\partial x^b}, \frac{\partial }{\partial x^a}),
\end{array},\eqno(6.14)$$
and {\bf R} is the scalar curvature of $(M, g)$. Then, we consider the {\it electric Weyl curvature tensor field}
$E = (E_{ac})$ given by

$$E_{ac} = E(\frac{\partial }{\partial x^c}, \frac{\partial }{\partial x^a}) = C_{abcd}\xi^b\xi^d,\eqno(6.15)$$
and taking into account that $\xi = \partial/\partial x^0$, we obtain

$$E_{ac} = C_{a0c0} = C\left(\frac{\partial }{\partial x^0}, \frac{\partial }{\partial
 x^c}, \frac{\partial }{\partial x^0}, \frac{\partial }{\partial x^a}\right).\eqno(6.16)$$
Then by direct calculations, using (6.16), (6.13), (2.7), (2.9) and (2.15), we deduce that the only possible
non-zero local components of $E$ with respect to the natural frame field are

$$E_{ik} = \bar{R}_{i0k0} + \frac{1}{2}\left\{\xi_i\bar{R}_{k0} + \xi_k\bar{R}_{i0}
- g_{ik}\bar{R}_{00} + \Phi^2\bar{R}_{ik}\right\} - \frac{1}{6}{\bf R}\Phi^2h_{ik}.\eqno(6.17)$$
Note that due to (2.3) and (6.16), we have

$$E_{ik} = C\left(\frac{\partial }{\partial x^0}, \frac{\delta }{\delta  x^k}, \frac{\partial }{\partial x^0},
 \frac{\delta }{\delta x^i}\right) = E\left(\frac{\delta }{\delta  x^k}, \frac{\delta }{\delta x^i}\right).
\eqno(6.18)$$
By using (2.5a) , from (6.18) we infer that $E_{ik}$ define a spatial tensor field of type
(0,2). Using (2.12) in (6.14b) and (6.14c), we obtain

$$\begin{array}{l}(a) \ \ \bar{R}_{i0k0} = R_{i0ko}, \ \ (b) \ \ \bar{R}_{i0} = R_{i0} - \Phi^{-2}\xi_iR_{00},
 \ \ (c) \ \ \bar{R}_{00} = R_{00}, \vspace{2mm} \\ (d) \ \ \ \bar{R}_{ik} = R_{ik} - \Phi^{-2}\left\{\xi_iR_{k0} +
 \xi_kR_{i0} - \Phi^{-2}\xi_i\xi_kR_{00}\right\}, \end{array}\eqno(6.19)$$
where $R_{i0k0}$ and $\{R_{ik}, R_{i0}, R_{00}\}$ are given by (5.2c) and (5.10), respectively. Taking into
account of (6.19) into (6.17) and using (2.15), we express $E_{ik}$ in terms of spatial tensor fields, as follows

$$E_{ik} = R_{i0k0} + \frac{1}{2}\left\{\Phi^2R_{ik} - \left(R_{00} +
\frac{1}{3}{\bf R}\Phi^2\right)h_{ik}\right\}.\eqno(6.20)$$
The scalar curvature ${\bf R}$ of $(M, g)$  is given by

$${\bf R} = \sum_{k=1}^3Ric(E_k, E_k) - \Phi^{-2}Ric\left(\frac{\partial }{\partial x^0},
\frac{\partial }{\partial x^0}\right) = h^{jh}R_{jh} - \Phi^{-2}R_{00}.\eqno(6.21)$$
Then, replace ${\bf R}$ from (6.21) into (6.20) and taking into account (6.11), we deduce that

$$E_{ik} = \widetilde{R}_{i0k0} + \frac{1}{2}\Phi^2\widetilde{R}_{ik}, \eqno(6.22)$$
where $\widetilde{R}_{ik}$ is the trace-free part of $R_{ik}$, that is, we have

$$\widetilde{R}_{ik} = R_{ik} - \frac{1}{3}h^{jh}R_{jh}h_{ik}. \eqno(6.23)$$
Finally, by using (6.11) and (6.22) into (6.8) and (6.10), we obtain

$$\begin{array}{l}
\Theta_{ik\vert_0} = \Phi^2\left\{b_ib_k + \frac{1}{2}\left(b_{i\vert_k} + b_{k\vert_i}\right)
 - \Phi^2\omega_{ih}\omega_k^h\right\} - \Theta _{ih}\Theta_k^h\vspace{2mm}\\ \hspace*{16mm}- E_{ik}
+ \frac{1}{2}\Phi^2\widetilde{R}_{ik} - \frac{1}{3} R_{00}h_{ik},\end{array}\eqno(6.24)$$
and

$$\begin{array}{l}
\sigma_{ik\vert_0} = \Phi^2\left\{b_ib_k + \frac{1}{2}\left(b_{i\vert_{k}} + b_{k\vert_{i}}\right)
- \frac{1}{3}\left(b^h_{\ \vert_h} + b^2\right)h_{ik}- \Phi^2\left(\omega_{ih}\omega_k^h\right.\right.
\vspace{2mm}\\ \left.\left.\hspace*{12mm}+
\frac{1}{3}\omega^2h_{ik}\right) \right\} + \frac{1}{3}\sigma^2h_{ik} - \sigma_{ih}\sigma^h_k - \frac{2}{3}
\Theta\sigma_{ik} - E_{ik} + \frac{1}{2}\Phi^2\widetilde{R}_{ik},\end{array}\eqno(6.25)$$
respectively.\par
It is interesting to note that the generalized Raychaudhuri's equation (6.1) can be expressed by using the
scalar curvature ${\bf R}$ of $(M, g)$ and the {\it spatial scalar curvature} ${\bf R^{\star}}$ of
$\nabla^{\star}$ given by

$${\bf R^{\star}} = h^{jh}R^{\star}_{jh}.\eqno(6.26)$$
Indeed, contracting (6.9) by $h^{ik}$ and using (3.16d), (3.8b), (6.21) and (6.26) we deduce that

$$\Theta_{\vert_0} = -\Theta^2 + \Phi^2\left\{b^2 + b^h_{\ \vert_h} + {\bf R} - {\bf R^{\star}}\right\}
 + R_{00}. \eqno(6.27)$$
In particular, if $\xi$ is a unit vector field that defines a timelike congruence of geodesics (see (6.3) and
(4.7)), then (6.27) becomes

$$\Theta_{\vert_0} = -\Theta^2 + {\bf R} - {\bf R^{\star}} + R_{00}. \eqno(6.28)$$
This is a {\it new form of the Raychaudhuri's equation} (6.6) {\it for a congruence of timelike geodesics}.
It is well known that (6.6) is the key equation used in the proof of Penrose-Hawking singularity theorems.
More precisely, it was proved the following Lemma.\par
{\bf Lemma 6.1.} (see Lemma 9.2.1 in \cite{w}) {\it Let $\xi$ be the tangent field of a hypersurface
orthogonal timelike geodesic congruence. Suppose the following conditions are satisfied:\par
(i) $Ric(\xi, \xi)\geq 0$, which is the case if Einstein's equations hold in the spacetime and the strong
energy condition is satisfied by the matter.\par
(ii) The expansion $\Theta$ takes the negative value $\Theta_0$ at a point on a geodesic in the
congruence corresponding to the proper time $\tau = 0.$\par
Then $\Theta$ goes to $-\infty$ along that geodesic within the proper time} $\tau \leq\frac{3}{\vert\Theta_0\vert}.$
Note that in the above lemma, $\xi$ is tangent field of a hypersurface orthogonal timelike geodesic congruence
means that $SM$ is an integrable distribution.\par
Now, by using the new form (6.28) of Raychaudhuri's equation we can complete Lemma 6.1 with the following Lemma.\par
{\bf Lemma 6.2.} {\it Let the congruence of timelike geodesics satisfying the conditions from Lemma 6.1. Then we
have the following assertions};\newline
(a) {\it If}  ${\bf R}\geq {\bf R^{\star}}$, {\it then the proper time $\tau$ must be in the interval}
 $[1/\vert\Theta_0\vert, 3/\vert\Theta_0\vert].$\newline
(b) \ \ {\it If} \ ${\bf R} < {\bf R^{\star}}$, {\it then the following cases occur}:\par
$(b_1)$ \ \  {\it If} \ $Ric(\xi, \xi) \geq {\bf R^{\star}} - {\bf R}$, {\it then $\tau$ must be in the interval}
 $[1/\vert\Theta_0\vert,  3/\vert\Theta_0\vert].$\par
$(b_2)$ \ \ {\it If} \ $Ric)\xi, \xi) < {\bf R^{\star}} - {\bf R}$, {\it then $\tau$ must be in the interval}
 $[0, 1/\vert\Theta_0\vert).$
{\bf Proof.} Suppose (a) is satisfied, and by using (i) in (6.28), we obtain

$$\Theta_{\vert_0} + \Theta^2 \geq 0,$$ which is equivalent to

$$\frac{d}{d\tau}(\frac{1}{\Theta})\leq 1.\eqno(6.29)$$
Then, integrating (6.29) on $[0, \tau]$, we infer that

$$\frac{1}{\Theta}\leq \frac{1}{\Theta_0} + \tau. \eqno(6.30)$$
As $1/\Theta$ must pass through zero, from (6.30) we deduce that
$$\tau \geq \frac{1}{\vert\Theta_0\vert}.\eqno(6.31)$$
Combining with the result from Lemma 6.1, we conclude that $\tau$ must be in the interval
$[1/\vert\Theta_0\vert,  3/\vert\Theta_0\vert].$ In a similar way it is proved the assertion $(b_1)$. Finally,
by using the condition from $(b_2)$  into (6.28), we obtain

$$\Theta_{\vert_0} + \Theta^2 < 0,$$ which is equivalent to

$$\frac{d}{d\tau}(\frac{1}{\Theta}) > 1.\eqno(6.32)$$
Integrating (6.32) on $[0, \tau]$, we deduce that

$$\frac{1}{\Theta}> \frac{1}{\Theta_0} + \tau. \eqno(6.33)$$
As $1/\Theta$ must pass through zero, we conclude that $\tau \in [0, 1/\vert\Theta_0\vert).$ This
completes the proof of the lemma. \ $\Box$\par

\section{Kinematic Quantities for Kerr-Newman Black Holes}

The new point of view developed here on the $(1+3)$ threading of spacetime, is applied in this section to the
charged Kerr black hole (also called Kerr-Newman black hole). We show that the curvature and Ricci tensor fields of
$(M,g)$ are simply expressed in terms of curvature and Ricci tensor fields of the Riemannian spatial connection,
via the kinematic quantities.\par
Now, according to the notations used in Sections 2 and 3, for the metric of a Kerr-Newman black hole given by
(1.2), we have

 $$\begin{array}{l}(a) \ \ \ \Phi^2 = \frac{\Delta - a^2 (\sin x^2)^2}{\Sigma} = 1 - \frac{2mx^1 - e^2}{\Sigma},
 \vspace{2mm}\\
 (b) \ \ \ \xi_1 = \xi_2 = 0, \ \ \xi_3 = \frac{(e^2 - 2mx^1)a(\sin x^2)^2}{\Sigma}, \vspace{2mm}\\
 (c) \ \ \ a_i = 0, \ \ \ \forall \ i\in \{1, 2, 3\}.
 \end{array}\eqno(7.1)$$
 The spatial distribution $SM$ of $(M, g)$ is locally spanned by

$$\frac{\delta }{\delta x^1} = \frac{\partial }{\partial x^1}, \ \ \frac{\delta }{\delta x^2} =
\frac{\partial }{\partial x^2}, \ \ \frac{\delta }{\delta x^3} = \frac{\partial }{\partial x^3} + \Phi^{-2} \xi_3
\frac{\partial }{\partial x^0},\eqno(7.2)$$
 and it is the kernel of the 1-form

$$\delta x^0 = dx^0 - \Phi^{-2}\xi_3dx^3. \eqno(7.3)$$
 By using (2.15) (1.2), (7.1a) and (7.1b), we deduce that the only non-zero local components
of the Riemannian metric $h$ on $SM$ with respect to the threading frame field from (7.2), are the following

$$h_{11} = \frac{\Sigma}{\Delta}, \ \ \ h_{22} = \Sigma, \ \ \ h_{33} = \frac{\Delta(\sin x^2)^2}{\Phi^2}.
\eqno(7.4) $$
Hence the line element from (1.2) becomes

$$ds^2 = -\Phi^2(\delta x^0)^2 + \frac{\Sigma}{\Delta}(dx^1)^2 + \Sigma(dx^2)^2 +
\frac{\Delta(\sin x^2)^2}{\Phi^2}(dx^3)^2,\eqno(7.5)$$
 with respect to the threading coframe field $\{\delta x^0, dx^i\}.$\par
 Next, by using (3.3), (4.6), (7.1a), and (7.1c), we deduce that the
geodesic spatial tensor field $b = (b_i)$ is given by

 $$\begin{array}{l}(a) \ \ \ b_1 = c_1 = \frac{x^1(2mx^1 - e^2) - m\Sigma}{(\Phi\Sigma)^2}, \vspace{2mm}\\
 (b) \ \ \ b_2 = c_2 = \frac{(e^2 - 2mx^1)a^2\sin x^2\cos x^2}{(\Phi\Sigma)^2}, \ \ \ (c) \ \ \ b_3 = c_3 = 0.
\end{array}\eqno(7.6)$$
 Due to (4.7) and (7.6) we conclude that the curves from the congruence defined by
$\xi = \partial/\partial x^0$ which sit in the surface given by the equations

$$ x^1(2mx^1 - e^2) - m\Sigma = 0, \ \ \  (e^2 - 2mx^1)a^2\sin x^2\cos x^2 = 0,$$
for $a\neq 0$, or in the hypersurface

$$x^1 = \frac{e^2}{m},$$
for $a = 0$, are the only geodesics of $(M, g)$ that are tangent to $\xi$. Moreover, we see that such geodesics
have the equations

$$\begin{array}{l} x^0 = \lambda, \ \ \ x^1 = \frac{e^2}{m}, \ \ \ x^2 = \frac{\pi}{2}, \ \ \ x^3 = c, \ \ \
\mbox{or}
\vspace{2mm}\\ x^0 = \lambda, \ \ \ x^1 = \frac{e^2 \pm \sqrt{e^4 + 4m^2}}{2m}, \ \ \ x^2 = 0, \ \ \ x^3 = c,
\end{array}$$
for $a \neq 0,$ and

 $$x^0 = \lambda, \ \ \ x^1 = \frac{e^2}{m}, \ \ \ x^2 = k, \ \ \ x^3 = c,$$
 for $a = 0$, where $k$ and $c$ are constants. In particular, the integral curves of $\xi$ can not be geodesics
in the Schwarzschild spacetime.\par
 Now, taking into account (3.8) and (7.4), we obtain

 $$(a) \ \ \ \Theta_{ij} = 0, \ \ \ (b) \ \ \ \Theta = 0, \ \ \ (c) \ \ \ \sigma_{ij} = 0, \ \ \
\forall \ i,j \in \{1, 2, 3\}.\eqno(7.7)$$
Also, by using (3.5a), (7.1b), (7.2) and (7.6), we deduce that the only non-zero local components of the vorticity tensor field $(\omega_{ij})$, are given by

$$\begin{array}{l}(a) \ \ \ \omega_{13} = \frac{a(\sin x^2)^2}{\Phi^4\Sigma^2}\{m\Sigma - x^1(2mx^1 - e^2)\},
\vspace{2mm}\\
(b) \ \ \ \omega_{23} = \frac{(2mx^1 - e^2)a\Delta\sin x^2\cos x^2}{\Phi^4\Sigma^2}.
 \end{array}\eqno(7.8)$$
 From (7.8) we see that the spatial distribution $SM$ is not integrable for both Kerr-Newman and Kerr black holes.
On the contrary, for the Reissner-Nordstrom and Schwarzschild solutions, the timelike vector field
$\xi = \partial/\partial x^0$ is hypersurface orthogonal.\par
Finally, we note that the local components of the curvature and Ricci tensor fields with respect to the
threading frame field have very simple expressions. Indeed, by using (7.6) and (7.7) into (5.3a),
(5.3b) and (5.5a), we obtain

$$\begin{array}{l}(a) \ \ \ R_{ijkh} = R^{\star}_{ijkh} + \Phi^2\left\{\omega_{ik}\omega_{jh} - \omega_{ih}\omega_{jk}
\right\},\vspace{2mm}\\
(b) \ \ \ R_{i0kh} = \Phi^2\left\{\omega_{ih\vert_k} - \omega_{ik\vert_h} + \omega_{ih}c_k - \omega_{ik}c_h +
2\omega_{kh}c_i\right\}, \vspace{2mm}\\
(c) \ \ \ R_{i0k0} = \Phi^2\left\{c_ic_k + \frac{1}{2}\left(c_{i\vert_k} + c_{k\vert_i}\right)
- \Phi^2\omega_{ih}\omega_k^h\right\},\end{array}\eqno(7.9)$$
where $R^{\star}_{ijkh}$ is the curvature tensor field of the Riemannian spatial connection. In a similar way,
from (5.12a), (5.11b) and (5.11c), we infer that the local components of the Ricci tensor of
$(M, g)$ with respect to the threading frame field  are given by

$$\begin{array}{l}(a) \ \ \ R_{ik} = R^{\star}_{ik} - c_ic_k - \frac{1}{2}\left(c_{i\vert_k} + c_{k\vert_i}\right),
\vspace{2mm}\\
(b) \ \ \ R_{i0} = \Phi^2\{\omega^h_{i\vert_h} + 3\omega_i^hc_h\},\vspace{2mm}\\
(c) \ \ \ R_{00} = \Phi^2\{c_hc^h + c^h_{\ \vert_h} + \Phi^2\omega_{kh}\omega^{kh}\},
\end{array}\eqno(7.10)$$
where $R^{\star}_{ik}$ is the Ricci tensor of the Riemannian spatial connection. As far as we know, (7.9) and (7.10)
have been not stated in earlier literature. They can bring more information and ideas in the study of the geometry
and physics
of the black holes. In particular, for Reissner-Nordstrom and Schwarzschild solutions, (7.9) and (7.10) become

$$\begin{array}{l}(a) \ \ \ R_{ijkh} = R^{\star}_{ijkh}, \ \ \ (b) \ \ \ R_{i0kh} = 0, \vspace{2mm} \\ (c) \ \ \
R_{i0k0} = \Phi^2\{c_ic_k + \frac{1}{2}\left(c_{i\vert_k} + c_{k\vert_i}\right)\},\end{array}\eqno(7.11)$$
and

$$\begin{array}{l}(a) \ \ \ R_{ik} = R_{ik}^{\star} - c_ic_k - \frac{1}{2}\left(c_{i\vert_k} + c_{k\vert_i}\right),
\vspace{2mm}\\
(b) \ \ \ R_{i0} = 0, \ \ \ (c) \ \ \ R_{00} = \Phi^2\{c_hc^h + c^h_{\ \vert_h}\}, \end{array}\eqno(7.12)$$
respectively.

\section{Equations of Motion in a Kerr Black Hole}

In this section and in the next one, we take $e = 0$ in (1.2), that is, we consider the Kerr black hole $(M,g)$.
As it is well known, the geodesic equations in a Kerr black hole have been explicitly integrated for first time by
Carter \cite{bm}. In the present section we will state a new form of the equations of motion in a Kerr black hole,
and obtain informations about the position of geodesics in $M$ with respect to the spatial distribution.
In particular, we show that the geodesics of $(M, g)$ which are tangent to $SM$ coincide with the autoparallel curves
of the Riemannian spatial connection.\par
Let $C$ be a smooth curve in $M$ given by parametric equations

$$x^0 = x^0(\lambda), \ \ \  x^i = x^i(\lambda), \ \ \ i \in \{1, 2, 3\}, \ \ \ \lambda \in [a,b], \eqno(8.1) $$
where $\lambda$ does not necessarily represents the time in $(M, g)$. The velocity vector field $d/d\lambda$ for
$C$ is expressed in terms of the threading frame $\{\partial/\partial x^0, \delta/\delta x^i\}$ as follows:

$$\frac{d }{d\lambda} = \frac{\delta x^0}{\delta\lambda}\frac{\partial }{\partial x^0}
+ \frac{dx^i}{d\lambda}\frac{\delta }{\delta x^i}, \eqno(8.2)$$
where we put

$$\frac{\delta x^0}{\delta\lambda} = \frac{dx^0}{d\lambda} - \Phi^{-2}\xi_3\frac{dx^3}{d\lambda}. \eqno(8.3)$$
Now, by using (7.7) and (7.1c) in (3.17), we express the Levi-Civita connection $\nabla$ on $(M, g)$,
as follows:

$$\begin{array}{lc}
(a) \ \ \nabla_{\frac{\delta }{\delta x^j}}\frac{\delta }{\delta x^i} =  \Gamma_{i \ j}^{\star k}
\frac{\delta }{\delta x^k} + \omega_{ij}\frac{\partial }{\partial x^0},
\vspace{4mm}\\
(b) \ \ \nabla_{\frac{\partial }{\partial x^0}}\frac{\delta }{\delta x^i} = \nabla_\frac{\delta }{\delta x^i}
{\frac{\partial }{\partial x^0}} = \Phi^2\omega_i^k\frac{\delta }{\delta x^k} +
c_i\frac{\partial }{\partial x^0}, \vspace{4mm}\\
(c) \ \ \nabla_{\frac{\partial }{\partial x^0}}\frac{\partial }{\partial x^0} =
 \Phi^2c^k\frac{\delta }{\delta x^k}. \end{array}\eqno(8.4)$$
Then, by direct calculations using (8.4) and (8.2), we obtain

$$\begin{array}{cr}\nabla_{\frac{d }{d\lambda}}\frac{d }{d\lambda}
= \left\{ \frac{d^2x^{k}}{d\lambda^2} + \Gamma_{i \ \ j}^{\star \ k}\frac{dx^{i}}{d\lambda}
 \frac{dx^{j}}{d\lambda} + 2\Phi^2\frac{\delta x^0}{\delta \lambda}\omega_{i}^{k}\frac{dx^{i}}{d\lambda}
+ \left(\frac{\delta x^0}{\delta \lambda}\right)^2\Phi^2c^{k}\right\}\frac{\delta }{\delta x^{k}}
\vspace{2mm}\\ + \left\{\frac{d }{d\lambda}(\frac{\delta x^0}{\delta \lambda}) +
2\frac{\delta x^0}{\delta \lambda}c_i\frac{dx^{i}}{d\lambda}\right\}\frac{\partial }{\partial x^0},\end{array}$$
which leads to the following equations of motion

$$\begin{array}{lcr} (a) \ \ \
\frac{d^2x^{k}}{d\lambda^2} + \Gamma_{i \ \ j}^{\star \ k}\frac{dx^{i}}{d\lambda}
 \frac{dx^{j}}{d\lambda} + 2\Phi^2\frac{\delta x^0}{\delta \lambda}\omega_{i}^{k}\frac{dx^{i}}{d\lambda}
+ \left(\frac{\delta x^0}{\delta \lambda}\right)^2\Phi^2c^{k} = 0,
\vspace{3mm}\\ (b) \ \ \ \frac{d }{d\lambda}(\frac{\delta x^0}{\delta \lambda}) +
2\frac{\delta x^0}{\delta \lambda}c_i\frac{dx^{i}}{d\lambda} = 0. \end{array}\eqno(8.5)$$
Note that $\lambda$ from (8.5) is an affine parameter for geodesics in $(M, g)$.\par
A geodesic of $(M, g)$ which is tangent at any of its points to the spatial distribution $SM$, is called a
{\it spatial geodesic}. Then by using (8.2), (8.3) and (8.5) we deduce that a curve $C$ given by (8.1) is a
spatial geodesic, if and only if, it is a solution of the system

$$\begin{array}{l}(a) \ \ \ \frac{d^2x^k}{d\lambda^2} + \Gamma_{i \ j}^{\star k}
 \frac{dx^i}{d\lambda}\frac{dx^j}{d\lambda} = 0,\vspace{3mm}\\
(b) \ \ \ \frac{\delta x^0}{\delta\lambda} = \frac{dx^0}{d\lambda} - \Phi^{-2}\xi_3\frac{dx^3}{d\lambda}
= 0.\end{array}\eqno(8.6)$$
Now, we remark that (7.7a) implies that the Kerr spacetime has bundle-like metric with respect to the
foliation determined by $\xi$ (cf.\cite{be}, p.112). Thus, we have the following interesting property:\par
{\it If a geodesic of a Kerr black hole is tangent to the spatial distribution at one point, then it remains
tangent to it at all later times}.\newline
Also, due to (8.6) we may state the following:\par
{\it The spatial geodesics in $(M, g)$ coincide with autoparallel curves for the Riemannian spatial connection}. \par
Next, we suppose that $C$ is a geodesic in $(M, g)$ which is not spatial, that is, we have
$\delta x^0/\delta\lambda \neq 0.$ Without loss of generality we suppose $\delta x^0/\delta\lambda > 0.$ Then by
using (3.3) in (8.5b), we obtain

$$\left(\frac{\delta x^0}{\delta\lambda}\right)^{-1}\frac{d }{d\lambda}\left(\frac{\delta x^0}{\delta\lambda}\right)
+ 2\Phi^{-1}\frac{d\Phi}{d\lambda} = 0, $$
which is equivalent to

$$\frac{\delta x^0}{\delta\lambda}\Phi^2 = K, \eqno(8.7)$$
where $K$ is a positive constant. By using (8.7) into (8.5a), we deduce that a geodesic of a Kerr black hole
(which is not a spatial geodesic), must be a solution of the system formed by (8.7) and the equations

$$\frac{d^2x^{k}}{d\lambda^2} + \Gamma_{i \ \ j}^{\star \ k}\frac{dx^{i}}{d\lambda}
 \frac{dx^{j}}{d\lambda} + 2K\omega_{i}^{k}\frac{dx^{i}}{d\lambda}
+ K^2\Phi^{-2}c^{k} = 0. \eqno(8.8)$$
Finally, taking into account (8.6b) and (8.7), we conclude that the system of differential equations for the
 geodesics in a Kerr black hole can be arranged in such a way that one of the equations is of first order.

\section{A $3D$ Identity Along a Geodesic in a Kerr Black Hole}

The Riemannian spatial connection given by (3.11), enables us to define a $3D$ force in a Kerr black hole,
and to deduce what we call the $3D$ force identity (cf.(9.10)). Note that this $3D$ force is a direct consequence
of the existence of the fourth dimension (time), and emphasizes an important difference between Newtonian gravity
and Einstein's general relativity.\par
Let $C$ be a geodesic in $(M, g)$, and $U(\lambda)$ be the projection of velocity $d/d\lambda$ on $SM$. Then by
(8.2) we obtain

$$U(\lambda) = \frac{dx^i}{d\lambda}\frac{\delta }{\delta x^i},\eqno(9.1) $$
which we call the $3D$ {\it velocity} along $C$. Now, consider the $3D$ arc-length parameter $s^{\star}$ on $C$
given by

$$s^{\star} = \int_a^{\lambda}h(U(\lambda), U(\lambda))^{1/2}d\lambda =
\int_a^{\lambda}\left(h_{ij}\frac{dx^i}{d\lambda}\frac{dx^j}{d\lambda}\right)^{1/2}d\lambda, $$
and obtain

$$(ds^{\star})^2 = h_{ij}dx^idx^j.\eqno(9.2)$$
Hence

$$U(s^{\star}) = \frac{dx^i}{ds^{\star}}\frac{\delta }{\delta x^i}, \eqno(9.3)$$
is a unit spatial vector field, that is, we have

$$h(U(s^{\star}), U(s^{\star})) = 1. \eqno(9.4)$$
Since $ds^{\star}/d\lambda$ is positive, we can take $s^{\star}$  as a new parameter on $C$. Then, define the $3D$
{\it force} along $C$, as the spatial vector field $F(s^{\star})$ given by

$$ F(s^{\star}) = \nabla^{\star}_{\frac{d }{ds^{\star}}}U(s^{\star}). \eqno(9.5)$$
Here, $\nabla^{\star}$ is the Riemannian spatial connection, and $d/ds^{\star}$ is given by

$$\frac{d }{ds^{\star}} = \frac{\delta x^0}{\delta s^{\star}}\frac{\partial }{\partial x^0} +
\frac{dx^i}{ds^{\star}}\frac{\delta }{\delta x^i}. \eqno(9.6)$$
Taking into account that $\nabla^{\star}$ is a metric connection on $SM$, and using (9.4) and (9.5), we deduce
that $F(s^{\star})$ is orthogonal to both $U(s^{\star})$ and $U(\lambda)$. By (9.5), (9.6) and (8.6), we see that
the $3D$ force vanishes along a spatial geodesic $C$, if and only if, $s^{\star}$ is an affine parameter on $C$.\par
Next, we put

$$ F(s^{\star}) = F^k(s^{\star})\frac{\delta }{\delta x^k}, $$
and by using (9.5), (9.6), (3.12), (3.14) and (7.7a), obtain

$$F^k(s^{\star}) = \frac{d^2x^{k}}{(ds^{\star})^2} + \Gamma_{i \ \ j}^{\star \ k}\frac{dx^{i}}{ds^{\star}}
 \frac{dx^{j}}{ds^{\star}}
+ \Phi^2\frac{\delta x^0}{\delta s^{\star}}\omega_{i}^{k}\frac{dx^i}{ds^{\star}}, \eqno(9.7)$$
provided $C$ is not a spatial geodesic. Now, by using  (8.7) and taking into account that

$$\frac{d^2\lambda}{(ds^{\star})^2} = -\frac{d^2s^{\star}}{d\lambda^2}\left(\frac{ds^{\star}}{d\lambda}\right)^{-3}, $$
from (9.7) we deduce that the local components of the $3D$ force with respect to the affine parameter
$\lambda$ , are given by

$$\begin{array}{r}
F^k(\lambda) = \left\{\frac{d^2x^{k}}{d\lambda^2} + \Gamma_{i \ \ j}^{\star \ k}\frac{dx^{i}}{d\lambda}
 \frac{dx^{j}}{d\lambda} + K\omega_{i}^{k}\frac{dx^{i}}{d\lambda}\right.\vspace{2mm}\\\left.
- \left(\frac{ds^{\star}}{d\lambda}\right)^{-1}
 \frac{d^2s^{\star}}{d\lambda^2}\frac{dx^k}{d\lambda}\right\}\left(\frac{ds^{\star}}{d\lambda}\right)^{-2}.
 \end{array}\eqno(9.8)$$
 Finally, by using (8.8) into (9.8), we infer that

$$F^k(\lambda) =  -\left\{K\omega_{i}^{k}\frac{dx^{i}}{d\lambda} + K^2\Phi^{-2}c^k
+ \left(\frac{ds^{\star}}{d\lambda}\right)^{-1}
 \frac{d^2s^{\star}}{d\lambda^2}\frac{dx^k}{d\lambda}\right\}\left(\frac{ds^{\star}}{d\lambda}\right)^{-2}. \eqno(9.9)$$
Taking into account that

$$h_{kh}F^k(\lambda)\frac{dx^h}{d\lambda} = 0,$$
and using (9.9),  (3.3) and (9.2), we obtain

$$K^2\Phi^{-3}\frac{d\Phi}{d\lambda} + \frac{ds^{\star}}{d\lambda}\frac{d^2s^{\star}}{d\lambda^2} = 0.\eqno(9.10)$$
Note that the identity (9.10) is a direct consequence of the existence of the $3D$ force $F$ given by (9.5). For
this reason we call it the $3D$ {\it force identity}. Such an identity could be useful in a study of motions in
$(M, g)$, and even for solving the equations of motion. For example, from (9.10) we deduce that:\par
 {\it The $3D$ arc length parameter $s^{\star}$ is an affine parameter on the geodesic $C$, if and only if, $\Phi$
is constant along $C$}.

\section{Conclusions}
In the present paper we develop a theory for a (1+3) threading of spacetime $(M,g)$ with respect to a congruence
of curves determined by an arbitrary timelike vector field $\xi = \partial/\partial x^0$. The generality of the
study is not the only difference between our approach and what is known in literature for the case of the unit
vector field $\xi$. The main differences consist in the following:\par
(i) Our approach is entirely developed with geometric objects expressed by their local components with respect
to the threading frames  $\{\partial/\partial x^0, \delta/\delta x^i\}$ and threading coframes
$\{\delta x^0, dx^i\}$.\par
(ii) The spatial distribution is not supposed to be necessarily integrable, and therefore this theory can be
easily applied to the study of any cosmological model with non-zero vorticity.\par
(iii) All the equations and results we state, are expressed in terms of the spatial tensor fields (cf.(3.1)) and
their spatial and time covariant derivatives (cf.(3.15)) induced by the Riemannian spatial connection given by
(3.12).\par
(iv) In spite of the numerous papers published on (1+3) threading of the spacetime (cf.[1,5]), the
generalized Raychaudhuri's equations (6.1), (6.5) and (6.27) are stated here for the first time in literature.\par
(v) The proof of Lemma 6.2 which completes the well known Lemma 6.1, is entirely based on a new form of
Raychaudhuri's equation for a congruence of timelike geodesics (cf.(6.28)).\par
(vi) It is the first time in literature when the spatial geodesics of a Kerr black hole are investigated (see (8.6)
and the assertions which follow it).\par
(vii) The $3D$ force (9.5) and the $3D$ force identity (9.10) are new objects in the general theory of Kerr back
holes, and illustrate the differences between Newtonian gravity and Einstein's general relativity.

\end{document}